\documentclass{elsart}
\usepackage{amssymb}
\usepackage{amsfonts}
\usepackage{amsmath}

\input{tcilatex}
\begin{document}

\begin{frontmatter}%

\title{Smoothing effect for the regularized Schrödinger equation with non controlled orbits}%

\author{Lassaad Aloui}%

\address{Département de Mathématiques, Faculté des Sciences de Bizerte, Tunisie\\ 
Email: lassaad.aloui@fsg.rnu.tn}%

\collab{Moez Khenissi}%

\address{Département de Mathématiques, École Supérieure des Sciences
et de Technologie de Hammam Sousse, Rue Lamine El Abbessi, 4011 Hammam Sousse, Tunisie\\ 
Email: moez.khenissi@fsg.rnu.tn}%

\collab{Georgi Vodev}%

\address{Université de Nantes, Département de Mathématiques, UMR 6629 du CNRS, 2, rue de la Houssinière, BP 92208, 44332 Nantes Cedex 03, France
\\ e-mail: georgi.vodev@math.univ-nantes.fr}%

\begin{abstract}
We prove that the geometric control condition is not necessary to obtain the
smoothing effect and the uniform stabilization for the strongly dissipative
Schr\"{o}dinger equation. 
\end{abstract}%


\end{frontmatter}%

\section{Introduction and statement of results}

It is well known that the Schr\"{o}dinger equation enjoys some smoothing
properties. One of them says that if $u_{0}\in L^{2}(\mathbb{R}^{d})$ with
compact support, then the solution of the Schr\"{o}dinger equation 
\begin{equation}
\left\{ 
\begin{array}{lll}
i\partial _{t}u-\Delta u=0 & \text{in} & \mathbb{R}\times \mathbb{R}^{d} \\ 
u(0,.)=u_{0} & \text{in} & \mathbb{R}^{d},%
\end{array}%
\right.   \label{nonreg}
\end{equation}%
satisfies 
\begin{equation*}
u\in C^{\infty }(\mathbb{R}\setminus \{0\}\times \mathbb{R}^{d}).
\end{equation*}%
We say that the Schr\"{o}dinger propagator has an infinite speed. Another
type of gain of regularity for system (\ref{nonreg}) is the Kato-$1/2$
smoothing effect (see \cite{co.sa1}, \cite{sjo}, \cite{veg}), namely any
solution of (\ref{nonreg}) satisfies 
\begin{equation}
\int_{\mathbb{R}}\int_{|x|<R}|(1-\Delta )^{\frac{1}{4}}u(t,x)|^{2}dxdt\leq
C_{R}\left\Vert u_{0}\right\Vert _{L^{2}(\mathbb{R}^{d})}^{2}.
\end{equation}%
In particular, this result implies that for a.e. $t\in \mathbb{R}$, $u(t,.)$
is locally smoother than $u_{0}$ and this happens despite the fact that (\ref%
{nonreg}) conserves the global $L^{2}$ norm. The Kato-effect has been
extended to variable coefficients operators with non trapping metric by Doi (%
\cite{doi1}, \cite{doi2})) and to non trapping exterior domains by Burq,
Gerard and Tzvetkov \cite{b.g.t}. On the other hand, Burq \cite{bursmot}
proved that the nontrapping assumption is necessary for the $H^{1/2}$
smoothing effect. Moreover, using Ikawa's result \cite{ik}, he showed, in
the case of several convex bodies satisfying certain assumptions, that the
smoothing effect with an $\varepsilon >0$ loss still holds.

Recently, the first author \cite{al1,al2} has introduced the forced
smoothing effect for Schr\"{o}dinger equation. The idea is inspired from the
stabilization problem and it consists of acting on the equation in order to
produce some smoothing effects. More precisely, in \cite{al2} the following
regularized Schr\"{o}dinger equation on a bounded domain $\Omega \subset 
\mathbb{R}^{d}$ is considered: 
\begin{equation}
\left\{ 
\begin{array}{lll}
i\partial _{t}u-\Delta _{D}u+ia(x)(-\Delta _{D})^{\frac{1}{2}}a(x)u=0 & 
\text{in} & \mathbb{R}\times \Omega , \\ 
u(0,.)=f & \text{in} & \Omega , \\ 
u|_{\mathbb{R}\times \partial \Omega }=0, &  & 
\end{array}%
\right.   \label{eqr}
\end{equation}%
where $\Delta _{D}$ denotes the Dirichlet realization of the Laplace
operator on $\Omega $ and $a(x)$ is a smooth real-valued function. Under the
geometric control condition (GCC) on the set $w=\{a\neq 0\}$, it is proved
in \cite{al2} that any solution with initial data in $H_{D}^{s}(\Omega )$
belongs to $L_{loc}^{2}((0,\infty ),H_{D}^{s+1}(\Omega ))$, where $s\in
\lbrack -s_{0},s_{0}]$ and $s_{0}\geq 1$ depends on the behavior of $a(x)$
near the boundary. When the function $a$ is constant near each component of
the boundary, we have $s_{0}=\infty $. Then by iteration of the last result,
a $C^{\infty }$-smoothing effect is proved in \cite{al2}. Note that these
smoothing effects hold away from $t=0$ and they seem strong compared with
the Kato effect for which the GCC is necessary. Therefore the case when $%
w=\{a\neq 0\}$ does not control geometrically $\Omega $ is very interesting.

In this work we give an example of geometry where the geometric control
condition is not satisfied but the $C^{\infty }$ smoothing effect holds.
More precisely, let $O=\cup _{i=1}^{N}O_{i}\subset \mathbb{R}^{d}$ be the
union of a finite number of bounded strictly convex bodies, $O_{i}$,
satisfying the conditions of \cite{ik}, namely:

\begin{itemize}
\item For any $1\leq i,j,k\leq N$, $i\neq j$, $j\neq k$, $k\neq i$, one has 
\begin{equation}
\text{Convex Hull}(O_{i}\cup O_{j})\cap O_{k}=\emptyset .
\end{equation}

\item Denote by $\kappa$ the infimum of the principal curvatures of the
boundaries of the bodies $O_{i}$, and $L$ the infimum of the distances
between two bodies. Then if $N>2$ we assume that $\kappa L>N$ (no assumption
if $N=2$).
\end{itemize}

Let $B$ be a bounded domain containing $O$ with smooth boundary and such
that $\Omega _{0}=O^{c}\cap B$ is connected, where $O^{c}=\mathbb{R}%
^d\setminus O$. In the present paper we will consider the regularized Schr%
\"{o}dinger equation (\ref{eqr}) in $\Omega _{0}$. For a bounded domain $%
\Omega $ of $\mathbb{R}^{d}$ and any $s\in \mathbb{R}$, we denote by $%
H_{D}^{s}(\Omega )$ the Hilbert space 
\begin{equation*}
H_{D}^{s}(\Omega )=\{u=\sum_{j}a_{j}e_{j},\quad \sum_{j}\gamma
_{j}^{2s}|a_{j}|^{2}<\infty \},
\end{equation*}
where $\{\gamma _{j}^{2}\}$ are the eigenvalues of $-\Delta _{D}$ and $%
\{e_{j}\}$ is the corresponding orthonormal basis of $L^{2}(\Omega )$. We
have the following interpolation inequalities:%
\begin{equation}
\text{ }\Vert g\Vert _{H_{D}^{s}(\Omega _{0})}\leq \Vert g\Vert
_{H_{D}^{t}(\Omega _{0})}^{\frac{s}{t}}\Vert g\Vert _{L^{2}(\Omega _{0})}^{1-%
\frac{s}{t}}\text{ for all }g\in H_{D}^{t}(\Omega _{0}),\quad 0\leq s\leq t.
\label{iterpol}
\end{equation}%
Clearly, $H_{D}^{s}(\Omega )$ and $H_{D}^{-s}(\Omega )$ are in duality and $%
H_{D}^{s}(\Omega )$ is the domain of $(-\Delta _{D})^{\frac{s}{2}}$. Remark
also that $H_{D}^{s}(\Omega )=H^{s}(\Omega )$ is the usual Sobolev space for 
$0\leq s<\displaystyle\frac{1}{2}$ and $H_{D}^{s}(\Omega )=\displaystyle%
\{u\in H^{s}(\Omega ),$ $\Delta ^{j}u|_{\partial \Omega }=0, 2j\leq s-\frac{1%
}{2}\}$ for $s\geq \displaystyle\frac{1}{2}$. Throughout this paper $a\in
C^\infty(\Omega_0)$ will be a real-valued function such that supp$%
\,a\subset\{x\in\overline{\Omega}_0:\mathrm{dist(x,\partial B)\le
2\varepsilon_0\}}$ and $a=Const\neq 0$ on $\{x\in\overline{\Omega}_0:\mathrm{%
dist(x,\partial B)\le \varepsilon_0\}}$, where $0<\varepsilon_0\ll 1$ is a
constant. Under this assumption the following properties hold for all $s\in 
\mathbb{R}$ and $n\in \mathbb{N}$:

\begin{description}
\item[ $\mathbf{(\mathcal{P}_{s})}$] the multiplication by $a$ maps $%
H^{s}_{D}(\Omega_{0})$ into itself,

\item[ $\mathbf{(\mathcal{Q}_{s,n})}$] the commutator $[a,(-\Delta_{D})^{n}]$
maps $H^{s}_{D}(\Omega_{0})$ into $H^{s-2n+1}_{D}(\Omega_{0})$.
\end{description}

Set $B_{a}=a(x)(-\Delta _{D})^{\frac{1}{2}}a(x)$ and define the operator $%
A_{a}=-\Delta _{D}+iB_a$ on $L^{2}(\Omega _{0})$ with domain 
\begin{equation*}
D(A_{a})=\{f\in L^{2}(\Omega _{0});A_{a}f\in L^{2}(\Omega _{0}),\quad
f=0\quad \mathrm{on}\quad \partial \Omega _{0}\}.
\end{equation*}%
Since the properties $\mathbf{(\mathcal{P}_{s})}$ and $\mathbf{(\mathcal{Q}
_{s,n})}$ hold for all $s\in \mathbb{R}$ and $n\in \mathbb{N}$, the problem
( \ref{eqr}) is well posed in $H_{D}^{s}(\Omega _{0})$ for all $s\in \mathbb{%
R} $. Moreover the operator $A_{a}$ generates a semi-group, $U(t)$, such
that for $f\in H_{D}^{s}(\Omega _{0}),$ $U(t)f\in C([0,+\infty \lbrack
,H_{D}^{s}(\Omega _{0}))$ is the unique solution of (\ref{eqr}). It is easy
to see that the spectrum, $\mathrm{sp}(A_{a})$, of $A_{a}$ consists of
complex numbers, $\tau _{j}$, satisfying $|\tau _{j}|\rightarrow \infty $.
Furthermore, since $a(x)$ is not identically zero, we have 
\begin{equation*}
\mathrm{sp}(A_{a})\subset \{\tau \in \mathbb{C},\text{Im}\,\tau >0\}.
\end{equation*}%
The resolvent 
\begin{equation*}
(A_a-\tau)^{-1}: L^{2}(\Omega _{0})\to L^{2}(\Omega _{0})
\end{equation*}
is holomorphic on $\ \left\{ \mathrm{Im}\,\tau <0\right\} $ and can be
extended to a meromorphic operator on $\mathbb{C}$. Our main result is the
following

\begin{theorem}
\label{Mainresult1} If the function $a$ is as above, there exist positive
constants $\sigma _{0}$ and $C$ such that for $|\mathrm{Im}\,\tau |<\sigma
_{0}$ we have 
\begin{equation}
\left\Vert (A_a-\tau )^{-1}\right\Vert _{L^{2}(\Omega _{0})\to L^{2}(\Omega
_{0})}\leq C\langle \tau \rangle ^{ -\frac{1}{2}}\log ^{2}\langle \tau
\rangle,  \label{res.hf}
\end{equation}
where $\langle \tau \rangle =\sqrt{1+|\tau |^{2}}$.
\end{theorem}

A similar bound has been recently proved in \cite{ca.vo} for the Laplace
operator in $\Omega_0$ with strong dissipative boundary conditions on $%
\partial B$, provided $B$ is strictly convex (viewed from the exterior).
Note also that a better bound (with $\log$ instead of $\log^2$) was obtained
in \cite{ch1}, \cite{ch2} in the case of the damped wave equation on compact
manifolds without boundary under the assumption that there is only one
closed hyperbolic orbit which does not pass through the support of the
dissipative term. This has been recently improved in \cite{sc} for a class
of compact manifolds with negative curvature, where a strip free of
eigenvalues has been obtained under a pressure condition.

As an application of this resolvent estimate we obtain the following
smoothing result for the associated Schr\"{o}dinger propagator.

\begin{theorem}
\label{Mainresult2} Let $s\in \mathbb{R}$. Under the assumptions of Theorem %
\ref{Mainresult1}, we have

(i) For each $\varepsilon >0$ there is a constant $C>0$ such that the
function 
\begin{equation*}
u(t)=\displaystyle\int^{t}_{0}e^{i(t-\tau)A_{a}} f(\tau)d\tau
\end{equation*}
satisfies 
\begin{equation}
\left\| u\right\| _{L^{2}_{T} H_{D}^{s+1-\varepsilon}(\Omega_{0} )}\leq
C\left\|f\right\| _{L^{2}_{T} H^{s}_{D}(\Omega_{0} )}  \label{regs}
\end{equation}
for all $T>0$ and $f\in L^{2}_{T}H^{s}_{D}(\Omega_{0} )$.

(ii)If $v_{0}\in H^{s}_{D}(\Omega_{0} )$, then 
\begin{equation}
v\in C^{\infty}((0,+\infty)\times \Omega_{0})  \label{regd}
\end{equation}
where $v$ is the solution of (\ref{eqr}) with initial data $v_{0}$.
\end{theorem}

Theorem \ref{Mainresult1} also implies the following stabilization result.

\begin{theorem}
\label{Mainresult3} Under the assumptions of Theorem \ref{Mainresult1},
there exist $\alpha ,c>0$ such that for the solution $u$ of (\ref{eqr}) with
initial data $u_{0}$ in $L^{2}(\Omega _{0})$, we have 
\begin{equation}
\Vert u\Vert _{L^{2}(\Omega _{0})}\leq ce^{-\alpha t}\Vert u_{0}\Vert
_{L^{2}(\Omega _{0})},\text{ }\forall \text{\ }t>1.  \notag
\label{unifDecay}
\end{equation}
\end{theorem}

This result shows that we can stabilize the Schr\"{o}dinger equation by a
(strongly) dissipative term that does not satisfy the geometric control
condition of \cite{Leb}. In fact, to have the exponential decay above it
suffices to have the estimate (\ref{res.hf}) with a constant in the
right-hand side.

\section{Resolvent estimates}

This section is devoted to the proof of Theorem \ref{Mainresult1}. Since the
resolvent $(A_{a}-\tau )^{-1}$ is meromorphic on $\mathbb{C}$ and has no
poles on the real axes, it suffices to prove (\ref{res.hf}) for $|\tau |\gg 1
$. Let $u$ be a solution of the following equation 
\begin{equation}
\left\{ 
\begin{array}{lll}
(\Delta _{D}+\lambda ^{2}-ia(x)(-\Delta _{D})^{\frac{1}{2}}a(x))u=v & \text{%
in} & \Omega _{0}, \\ 
u|_{\partial \Omega _{0}}=0, &  & 
\end{array}%
\right.   \label{eqsemiclassique}
\end{equation}%
with $v\in L^{2}(\Omega _{0})$. Clearly, it suffices to prove the following

\begin{proposition}
\label{SemEst} Under the assumptions of Theorem \ref{Mainresult1}, there
exists $\lambda _{0}>0$ such that for every $\lambda >\lambda _{0}$ and
every solution $u$ of (\ref{eqsemiclassique}) we have 
\begin{equation}
\Vert u\Vert _{L^{2}(\Omega _{0})}\lesssim \frac{\log^2 \lambda}{\lambda }%
\Vert v\Vert _{L^{2}(\Omega _{0})}.
\end{equation}
\end{proposition}

\textit{Proof.} Let $\chi \in C_{0}^{\infty }(B)$ be such that $\chi=0 $ on $%
\{x\in B: \mathrm{dist}(x,\partial B)\le\varepsilon_0/3\}$, $\chi=1$ on $%
\{x\in B: \mathrm{dist}(x,\partial B)\ge\varepsilon_0/2\}$. Clearly we have 
\begin{equation}
\Vert (1-\chi )u\Vert _{L^{2}(\Omega _{0})}\lesssim \Vert au\Vert
_{L^{2}(\Omega _{0})}.  \label{a1}
\end{equation}
On the other hand, the function $\chi u$ satisfies the equation 
\begin{equation}
\left\{ 
\begin{array}{lll}
(\Delta _{D}+\lambda ^{2})\chi u=\chi v+i\chi a(x)(-\Delta _{D})^{\frac{1}{2}
}a(x)u+[\Delta _{D},\chi ]u & \quad\text{in} & \mathbb{R}^{d}\setminus O, \\ 
u|_{\partial O}=0. &  & 
\end{array}
\right.  \label{eqsemiext}
\end{equation}
Hence, according to Proposition 4.8 of \cite{bursmot}, it follows that 
\begin{equation}
\Vert \chi u\Vert _{L^{2}(\Omega _{0})}\lesssim \frac{\log \lambda }{\lambda 
}(\Vert v\Vert _{L^{2}(\Omega _{0})}+\Vert au\Vert _{H^{1}(\Omega _{0})}).
\label{a2}
\end{equation}
By (\ref{a1}) and (\ref{a2}), 
\begin{equation}
\Vert u\Vert _{L^{2}(\Omega _{0})}\lesssim \frac{\log \lambda }{\lambda }
\Vert v\Vert _{L^{2}(\Omega _{0})}+\frac{\log \lambda }{\lambda }\Vert
au\Vert _{H^{1}(\Omega _{0})}+\Vert au\Vert _{L^{2}(\Omega _{0})}.
\label{estLog}
\end{equation}
To estimate the second and the third terms in the right hand side of (\ref%
{estLog}), we need the following

\begin{lemma}
\label{lemHs}Let $s\in \lbrack 0,1]$ and $\psi \in C^{\infty }(\Omega _{0})$%
. Then we have, for $\lambda\gg 1$, 
\begin{equation}
\Vert \psi u\Vert _{H^{s+1}}\lesssim \lambda \Vert \psi u\Vert
_{H^{s}}+\Vert v\Vert _{L^{2}}+\lambda ^{1/2}\Vert u\Vert _{L^{2}},
\label{hs+1}
\end{equation}
\begin{equation}
\Vert \psi u\Vert _{H^{s}}\lesssim \lambda ^{-1}\Vert \psi u\Vert
_{H^{s+1}}+\lambda ^{-1}\Vert v\Vert _{L^{2}}+\lambda ^{-1/2}\Vert u\Vert
_{L^{2}}.  \label{hs}
\end{equation}
\end{lemma}

\textit{Proof.} The function $w=(-\Delta _{D})^{s/2}\psi u:=P_{s}u$
satisfies the equation 
\begin{equation}
\left\{ 
\begin{array}{lll}
(-\Delta _{D}-\lambda ^{2}+iB_{a})w=P_{s}v+[\Delta
_{D},P_{s}]u+i[B_{a},P_{s}]u & \text{in} & \Omega _{0}, \\ 
w|_{\partial O}=0. &  & 
\end{array}%
\right.  \label{Ps}
\end{equation}%
Multiplying equation (\ref{Ps}) by $\overline{w}$, integrating by parts and
taking the real part, we obtain 
\begin{equation}
\Vert (-\Delta _{D})^{1/2}w\Vert _{L^{2}(\Omega _{0})}^{2}-\lambda ^{2}\Vert
w\Vert _{L^{2}(\Omega _{0})}^{2}=\mathrm{Re}\,\langle P_{s}v+[\Delta
_{D},P_{s}]u+i[B_{a},P_{s}]u,w\rangle .  \label{eqW}
\end{equation}%
Using that $[\Delta _{D},P_{s}]=(-\Delta _{D})^{s/2}[\Delta _{D},\psi ]$, we
deduce from (\ref{eqW}) 
\begin{equation}
\begin{array}{l}
\Vert (-\Delta _{D})^{1/2}w\Vert _{L^{2}}^{2}\lesssim \lambda ^{2}\Vert
w\Vert _{L^{2}}^{2}+\Vert \psi v\Vert _{L^{2}}\Vert P_{2s}u\Vert
_{L^{2}}+\Vert u\Vert _{H^{s+1}}\Vert w\Vert _{L^{2}}.%
\end{array}%
\end{equation}%
This implies 
\begin{equation}
\begin{array}{l}
\Vert \psi u\Vert _{H^{s+1}}^{2} \\ 
\lesssim \lambda ^{2}\Vert \psi u\Vert _{H^{s}}^{2}+\Vert v\Vert
_{L^{2}}\Vert \psi u\Vert _{H^{2s}}+\Vert u\Vert _{H^{s+1}}\Vert \psi u\Vert
_{H^{s}} \\ 
\lesssim \lambda ^{2}\Vert \psi u\Vert _{H^{s}}^{2}+\Vert v\Vert
_{L^{2}}^{2}+\varepsilon \Vert \psi u\Vert _{H^{2s}}^{2}+\Vert u\Vert
_{H^{s+1}}\Vert \psi u\Vert _{H^{s}} \\ 
\lesssim \lambda ^{2}\Vert \psi u\Vert _{H^{s}}^{2}+\Vert v\Vert
_{L^{2}}^{2}+\varepsilon \Vert \psi u\Vert _{H^{s+1}}^{2}+\varepsilon \Vert
\psi u\Vert _{H^{s}}^{2}+\lambda ^{-2}\Vert u\Vert _{H^{s+1}}^{2}+\lambda
^{2}\Vert \psi u\Vert _{H^{s}}^{2}.%
\end{array}%
\end{equation}%
Choosing $\varepsilon $ small enough, we obtain 
\begin{equation}
\Vert \psi u\Vert _{H^{s+1}}^{2}\lesssim \lambda ^{2}\Vert \psi u\Vert
_{H^{s}}^{2}+\Vert v\Vert _{L^{2}}^{2}+\lambda ^{-2}\Vert u\Vert
_{H^{s+1}}^{2}.  \label{s-1}
\end{equation}%
Taking $\psi =1$ in the previous estimate, we get for $\lambda \gg 1$ 
\begin{equation}
\Vert u\Vert _{H^{s+1}}^{2}\lesssim \lambda ^{2}\Vert u\Vert
_{H^{s}}^{2}+\Vert v\Vert _{L^{2}}^{2}.  \label{psi1}
\end{equation}%
Inserting (\ref{psi1}) in (\ref{s-1}), we obtain 
\begin{equation}
\Vert \psi u\Vert _{H^{s+1}}^{2}\lesssim \lambda ^{2}\Vert \psi u\Vert
_{H^{s}}^{2}+\Vert v\Vert _{L^{2}}^{2}+\Vert u\Vert _{H^{s}}^{2}.
\label{psis}
\end{equation}%
On the other hand, since $s\in \lbrack 0,1]$, by interpolation we have 
\begin{equation}
\Vert u\Vert _{H^{s}}^{2}\leq \Vert u\Vert _{H^{1}}\Vert u\Vert _{L^{2}}\leq
\lambda ^{-1}\Vert u\Vert _{H^{1}}^{2}+\lambda \Vert u\Vert _{L^{2}}^{2}.
\label{h2}
\end{equation}%
Choosing $s=0$ in (\ref{psi1}), we get 
\begin{equation}
\Vert u\Vert _{H^{1}}^{2}\lesssim \lambda ^{2}\Vert u\Vert
_{L^{2}}^{2}+\Vert v\Vert _{L^{2}}^{2}.  \label{H1}
\end{equation}%
Combining (\ref{h2}) and (\ref{H1}), we obtain 
\begin{equation}
\Vert u\Vert _{H^{s}}^{2}\lesssim \lambda \Vert u\Vert _{L^{2}}^{2}+\lambda
^{-1}\Vert v\Vert _{L^{2}}^{2}.
\end{equation}%
Inserting this estimate in (\ref{psis}) we get 
\begin{equation}
\Vert \psi u\Vert _{H^{s+1}}^{2}\lesssim \lambda ^{2}\Vert \psi u\Vert
_{H^{s}}^{2}+\Vert v\Vert _{L^{2}}^{2}+\lambda \Vert u\Vert _{L^{2}}^{2}.
\end{equation}%
This completes the proof of (\ref{hs+1}). Clearly, (\ref{hs}) can be proved
in the same way.

We deduce from Lemma \ref{lemHs} the following

\begin{lemma}
\label{lempsi}Let $u$ be a solution of (\ref{eqsemiclassique}) and $\psi \in
C^{\infty }(\Omega _{0})$. Then we have, for $\lambda \gg 1$, 
\begin{equation}
\Vert \psi u\Vert _{L^{2}}\lesssim \lambda ^{-1/2}\Vert \psi u\Vert
_{H^{1/2}}+\lambda ^{-1}\Vert v\Vert _{L^{2}}+\lambda ^{-1/2}\Vert u\Vert
_{L^{2}},  \label{cor1}
\end{equation}%
\begin{equation}
\Vert \psi u\Vert _{H^{1}}\lesssim \lambda ^{1/2}\Vert \psi u\Vert
_{H^{1/2}}+\lambda ^{-1/2}\Vert v\Vert _{L^{2}}+\Vert u\Vert _{L^{2}}.
\label{cor2}
\end{equation}
\end{lemma}

\textit{Proof.} Using Lemma \ref{lemHs} together with an interpolation
argument, we get 
\begin{equation}
\begin{array}{lll}
\Vert \psi u\Vert _{H^{1}}^{2} & \lesssim & \lambda ^{-1}\Vert \psi u\Vert
_{H^{3/2}}^{2}+\lambda \Vert \psi u\Vert _{H^{1/2}}^{2} \\ 
& \lesssim & \lambda ^{-1}(\lambda ^{2}\Vert \psi u\Vert
_{H^{1/2}}^{2}+\Vert v\Vert _{L^{2}}^{2}+\lambda \Vert u\Vert
_{L^{2}}^{2})+\lambda \Vert \psi u\Vert _{H^{1/2}}^{2} \\ 
& \lesssim & \lambda \Vert \psi u\Vert _{H^{1/2}}^{2}+\lambda ^{-1}\Vert
v\Vert _{L^{2}}^{2}+\Vert u\Vert _{L^{2}}^{2},%
\end{array}%
\end{equation}
which proves (\ref{cor2}). Using (\ref{eqW}) and (\ref{cor2}) we obtain 
\begin{equation}
\begin{array}{lll}
\Vert \psi u\Vert _{L^{2}}^{2} & \lesssim & \lambda ^{-2}\Vert \psi u\Vert
_{H^{1}}^{2}+\lambda ^{-2}\Vert v\Vert _{L^{2}}^{2}+\lambda ^{-1}\Vert
u\Vert _{L^{2}}^{2} \\ 
& \lesssim & \lambda ^{-1}\Vert \psi u\Vert _{H^{1/2}}^{2}+\lambda
^{-2}\Vert v\Vert _{L^{2}}^{2}+\lambda ^{-1}\Vert u\Vert _{L^{2}}^{2},%
\end{array}%
\end{equation}
which proves (\ref{cor1}).

We now return to the proof of Proposition \ref{SemEst}. Using Lemma \ref%
{lempsi} and the estimate (\ref{estLog}), we get 
\begin{equation}
\Vert u\Vert _{L^{2}}\lesssim \frac{\log\lambda }{\lambda }\Vert v\Vert
_{L^{2}}+\frac{\log\lambda }{\sqrt{\lambda }}\Vert au\Vert _{H^{1/2}}\text{
for }\lambda \gg 1.  \label{estlogA}
\end{equation}
Let's now estimate the $H^{1/2}$ term. Multiplying equation (\ref%
{eqsemiclassique}) by $\overline{u}$, integrating by parts and taking the
imaginary part, we obtain 
\begin{equation}
\Vert au\Vert _{H^{1/2}}^{2}=\langle (-\Delta _{D})^{\frac{1}{2}
}au,au\rangle =\mathrm{Im}\,\langle v,u\rangle \leq \Vert v\Vert
_{L^{2}}\Vert u\Vert _{L^{2}}.  \label{est}
\end{equation}
By (\ref{estlogA}) and (\ref{est}), we get 
\begin{equation}
\Vert u\Vert _{L^{2}}\lesssim \frac{\log\lambda }{\lambda }\Vert v\Vert
_{L^{2}}+\frac{\log\lambda }{\sqrt{\lambda }}\Vert v\Vert
_{L^{2}}^{1/2}\Vert u\Vert _{L^{2}}^{1/2}.
\end{equation}
This implies 
\begin{equation}
\Vert u\Vert _{L^{2}}\lesssim \frac{\log\lambda }{\lambda }\Vert v\Vert
_{L^{2}}+\frac{\log\lambda }{\sqrt{\lambda }}(\frac{\log\lambda }{\sqrt{%
\lambda }}\Vert v\Vert _{L^{2}}+\varepsilon \frac{\sqrt{\lambda }}{
\log\lambda }\Vert u\Vert _{L^{2}}),
\end{equation}
for any $\varepsilon >0$. Choosing $\varepsilon$ small enough, we get for
large $\lambda $ 
\begin{equation}
\Vert u\Vert _{L^{2}}\lesssim \frac{\log^2\lambda }{\lambda }\Vert v\Vert
_{L^{2}},
\end{equation}
which is the desired result.

\section{Smoothing effect}

We will first prove the following

\begin{proposition}
\label{Resolvent-Hs} If $a(x)$ is as in the introduction, then for every $%
s\in \mathbb{R}$, $\varepsilon>0$ there exist positive constants $C$ and $%
\sigma_{0}$ such that 
\begin{equation}
\|(A_a-\tau)^{-1}\| _{H_{D}^{s}(\Omega_{0})\rightarrow
H_{D}^{s+1-\varepsilon}(\Omega_{0})}\leq C  \label{resso}
\end{equation}
holds for $|\mathrm{Im}\,\tau|<\sigma_{0}$.
\end{proposition}

\textit{Proof.} Let $u$ and $f$ satisfy the equation 
\begin{equation}
\left\{ 
\begin{array}{lll}
(-\Delta _{D}-\tau +iB_{a})u=f & \text{in} & \Omega _{0}, \\ 
u=0 & \text{on} & \partial \Omega _{0}.%
\end{array}%
\right.  \label{las}
\end{equation}%
Let's see that the following estimate holds 
\begin{equation}
\Vert u\Vert _{H_{D}^{2}(\Omega _{0})}\leq C\langle \tau \rangle ^{\frac{1}{2%
}}\log ^{2}\langle \tau \rangle \Vert f\Vert _{L^{2}(\Omega _{0})},
\label{resalpha}
\end{equation}%
for $|\mathrm{Im}\,\tau |<\sigma _{0}$. Using Proposition 1 we get 
\begin{equation}
\begin{array}{lll}
\Vert u\Vert _{H_{D}^{2}(\Omega _{0})} & = & \Vert \Delta _{D}u\Vert
_{L^{2}(\Omega _{0})} \\ 
& = & \Vert \tau u-iB_{a}u+f\Vert _{L^{2}(\Omega _{0})} \\ 
& \leq & C\langle \tau \rangle ^{\frac{1}{2}}\log ^{2}\langle \tau \rangle
\Vert f\Vert _{L^{2}(\Omega _{0})}+C\Vert u\Vert _{H_{D}^{1}(\Omega _{0})}
\\ 
& \leq & C\langle \tau \rangle ^{\frac{1}{2}}\log ^{2}\langle \tau \rangle
\Vert f\Vert _{L^{2}(\Omega _{0})}+\varepsilon \Vert u\Vert
_{H_{D}^{2}(\Omega _{0})}+C_{\varepsilon }\Vert f\Vert _{L^{2}(\Omega _{0})}.%
\end{array}%
\end{equation}%
Choosing $\varepsilon $ small enough, we obtain (\ref{resalpha}). Using (\ref%
{resalpha}) and (\ref{iterpol}) with $s=1-\varepsilon $, $t=2$, we obtain 
\begin{equation}
\begin{array}{lll}
\Vert u\Vert _{H_{D}^{1-\varepsilon }(\Omega _{0})}^{2} & \leq & %
\displaystyle\Vert u\Vert _{H_{D}^{2}(\Omega _{0})}^{1-\varepsilon }\Vert
u\Vert _{L^{2}(\Omega _{0})}^{1+\varepsilon } \\ 
& \leq & C\displaystyle(\langle \tau \rangle ^{\frac{1}{2}}\log ^{2}\langle
\tau \rangle )^{1-\varepsilon }(\frac{\log ^{2}\langle \tau \rangle }{%
\langle \tau \rangle ^{\frac{1}{2}}})^{1+\varepsilon }\Vert f\Vert
_{L^{2}(\Omega _{0})}^{2} \\ 
& \leq & C\Vert f\Vert _{L^{2}(\Omega _{0})}^{2}.%
\end{array}%
\end{equation}%
Then we get (\ref{resso}) for $s=0$, i.e. 
\begin{equation}
\Vert (-\Delta _{D}-\tau +iB_{a})^{-1}\Vert _{L^{2}(\Omega _{0})\rightarrow
H_{D}^{1-\varepsilon }(\Omega _{0})}\leq C.  \label{h1-epsilon}
\end{equation}%
We will prove (\ref{resso}) for $s=2N$ with $N\in \mathbb{N}$, namely 
\begin{equation}
\Vert u\Vert _{H_{D}^{2N+1-\varepsilon }(\Omega _{0})}\lesssim \Vert f\Vert
_{H_{D}^{2N}(\Omega _{0})}.  \label{hhd}
\end{equation}%
Let $f\in H_{D}^{2N}(\Omega _{0})$ and let $u$ be the corresponding solution
of (\ref{las}). The function $(-\Delta _{D})^{N}u$ satisfies 
\begin{equation}
(\Delta _{D}+\tau -iB_{a})((-\Delta _{D})^{N}u)=(-\Delta
_{D})^{N}f-i[B_{a},(-\Delta _{D})^{N}]u.
\end{equation}%
Using (\ref{h1-epsilon}), we obtain 
\begin{equation}
\Vert (-\Delta _{D})^{N}u\Vert _{H^{\gamma }(\Omega _{0})}\lesssim \Vert
(-\Delta _{D})^{N}f\Vert _{L^{2}(\Omega _{0})}+\Vert \lbrack B_{a},(-\Delta
_{D})^{N}]u\Vert _{L^{2}(\Omega _{0})}
\end{equation}%
where $\gamma =1-\varepsilon $. Since 
\begin{equation}
\Vert u\Vert _{H_{D}^{2N+\gamma }(\Omega _{0})}=\Vert (-\Delta )^{N}u\Vert
_{H^{\gamma }(\Omega _{0})},
\end{equation}%
we obtain 
\begin{equation}
\Vert u\Vert _{H^{2N+\gamma }(\Omega _{0})}\lesssim \Vert (-\Delta
_{D})^{N}f\Vert _{L^{2}(\Omega _{0})}+\Vert \lbrack B_{a},(-\Delta
_{D})^{N}]u\Vert _{L^{2}(\Omega _{0})}.
\end{equation}%
On the other hand, we have 
\begin{equation}
\lbrack B_{a},(-\Delta _{D})^{N}]=a(-\Delta _{D})^{\frac{1}{2}}[a,(-\Delta
_{D})^{N}]+[a,(-\Delta _{D})^{N}](-\Delta _{D})^{\frac{1}{2}}a.
\end{equation}%
Using the properties $\mathbf{(\mathcal{P}_{s})}$ and $\mathbf{(\mathcal{Q}%
_{s,n})}$ we get 
\begin{equation}
\Vert \lbrack B_{a},(-\Delta _{D})^{N}]u\Vert _{L^{2}(\Omega _{0})}\lesssim
\Vert u\Vert _{H_{D}^{2N}(\Omega _{0})}.
\end{equation}%
Consequently 
\begin{equation}
\begin{array}{lll}
\Vert u\Vert _{H_{D}^{2N+\gamma }(\Omega _{0})} & \lesssim & \Vert f\Vert
_{H_{D}^{2N}(\Omega _{0})}+\Vert u\Vert _{H_{D}^{2N}(\Omega _{0})} \\ 
& \lesssim & \Vert f\Vert _{H_{D}^{2N}(\Omega _{0})}+\varepsilon \Vert
u\Vert _{H_{D}^{2N+\gamma }(\Omega _{0})}+C_{\varepsilon }\Vert u\Vert
_{H_{D}^{\gamma }(\Omega _{0})}.%
\end{array}%
\end{equation}%
Choosing $\varepsilon $ small enough and using (\ref{h1-epsilon}) we obtain 
\begin{equation}
\begin{array}{ccl}
\Vert u\Vert _{H_{D}^{2N+\gamma }(\Omega _{0})} & \lesssim & \Vert f\Vert
_{H_{D}^{2N}(\Omega _{0})}+\Vert u\Vert _{H_{D}^{\gamma }(\Omega _{0})} \\ 
& \lesssim & \Vert f\Vert _{H_{D}^{2N}(\Omega _{0})}+\Vert f\Vert
_{L^{2}(\Omega _{0})} \\ 
& \lesssim & \Vert f\Vert _{H_{D}^{2N}(\Omega _{0})},%
\end{array}%
\end{equation}%
which proves (\ref{resso}) for $s=2N$. Using the identity 
\begin{equation}
\lbrack a,(-\Delta _{D})^{-N}]=(-\Delta _{D})^{-N}[(-\Delta
_{D})^{N},a](-\Delta _{D})^{-N},
\end{equation}%
we can prove (\ref{resso}) for $s=-2N$ in the same way as in the case $s=2N$%
. Finally, by an interpolation argument we obtain the result for $s\in 
\mathbb{R}$. This completes the proof of Proposition \ref{Resolvent-Hs}.

\textit{Proof of Theorem \ref{Mainresult2}.} We will first prove (\ref{regs}%
). Extend $f$ by $0$ for $t\in \mathbb{R}\setminus[0,T]$. The Fourier
transforms (in $t$) of $u$ and $f$ are holomorphic in the domain $\mathrm{Im}%
\,z<0$ and satisfy the equation 
\begin{equation}
(-z-\Delta_{D}+iB_{a})\widehat{u}(z,\cdot)=\widehat{f}(z,\cdot).
\end{equation}
We take $z=\lambda-i\sigma$, $\lambda\in\mathbb{R}$ , $\sigma>0$, and we let 
$\sigma$ tend to zero. Using Proposition \ref{Resolvent-Hs}, we get 
\begin{equation}
\| \widehat{u}\|_{L^{2}( \mathbb{R};
H^{s+1-\varepsilon}_{D}(\Omega_{0}))}\lesssim \|\widehat{f}\|_{L^{2}( 
\mathbb{R}; H^{s}_{D}(\Omega_{0}))},\qquad s\in \mathbb{R}.
\end{equation}
The fact that the Fourier transform of any function from $\mathbb{R}$ to a
Hilbert space $H$ defines an isometry on $L^{2}( \mathbb{R}; H)$ completes
the proof of (\ref{regs}).

Now we turn to the proof of (\ref{regd}). Let $\varphi \in C_{0}^{\infty
}]0,+\infty \lbrack ,$ then the function $w(t,\cdot )=\varphi (t)v(t,\cdot )$
satisfies the equation 
\begin{equation}
\left\{ 
\begin{array}{lll}
i\partial _{t}w-\Delta _{D}w+iB_{a}w=i\varphi ^{\prime }(t)v & \text{in} & 
\mathbb{R}_{+}\times \Omega _{0}, \\ 
w(0,\cdot )=0 & \text{in} & \Omega _{0}, \\ 
w|_{\mathbb{R}\times \partial \Omega _{0}}=0. &  & 
\end{array}%
\right.
\end{equation}%
Using (\ref{regs}) with $\varepsilon =1/2$ we obtain 
\begin{equation}
\begin{array}{lll}
\left\Vert w\right\Vert _{L^{2}(\mathbb{R}_{+},H_{D}^{s+1/2}(\Omega _{0}))}
& \lesssim & \left\Vert \varphi ^{\prime }v\right\Vert _{L^{2}(\mathbb{R}%
_{+},H_{D}^{s}(\Omega _{0}))} \\ 
& \lesssim & \left\Vert v_{0}\right\Vert _{H_{D}^{s}(\Omega _{0})},%
\end{array}%
\end{equation}%
which implies 
\begin{equation}
v\in L_{loc}^{2}((0,\infty ),H_{D}^{s+1/2}(\Omega _{0})).
\end{equation}%
By iteration we obtain 
\begin{equation}
v\in L_{loc}^{2}((0,\infty ),H_{D}^{s+k}(\Omega _{0})),\qquad \forall k\in 
\mathbb{N}.
\end{equation}%
Using the equation satisfied by $v,$ we deduce that 
\begin{equation}
v\in H_{loc}^{k}((0,\infty ),H_{D}^{s+k}(\Omega _{0})),\qquad \forall k\in 
\mathbb{N}.
\end{equation}%
This implies that $v\in C^{\infty }((0,\infty )\times \Omega _{0})$ and the
proof of Theorem \ref{Mainresult2} is completed.

\end{document}